\documentclass{article}%
\usepackage{amsfonts}
\usepackage{makeidx}
\usepackage{amsmath}
\usepackage{amssymb}
\usepackage{graphicx}%
\providecommand{\U}[1]{\protect\rule{.1in}{.1in}}
\newtheorem{theorem}{Theorem}

\begin{document}

\title{Lie and Noether point Symmetries for a Class of Nonautonomous Dynamical Systems}
\author{Leonidas Karpathopoulos\\{\ \ \textit{Faculty of Physics, Department of
Astronomy-Astrophysics-Mechanics,}}\\{\ \textit{University of Athens, Panepistemiopolis, Athens 157 83, Greece}}
\and Andronikos Paliathanasis\thanks{Email: anpaliat@phys.uoa.gr}\\{\ \textit{Instituto de Ciencias F\'{\i}sicas y Matem\'{a}ticas, }}\\{\textit{Universidad Austral de Chile, Valdivia, Chile}}\\{\textit{Institute of Systems Science, Durban University of Technology}}\\{\textit{Durban 4000, Republic of South Africa}}
\and Michael Tsamparlis\thanks{Email: mtsampa@phys.uoa.gr}\\{\ \textit{Faculty of Physics, Department of
Astronomy-Astrophysics-Mechanics,}}\\{\ \textit{University of Athens, Panepistemiopolis, Athens 157 83, Greece}}}
\date{}
\maketitle

\begin{abstract}
We prove two general theorems which determine the Lie and the Noether point
symmetries for the equations of motion of a dynamical system which moves in a
general Riemannian space under the action of a time dependent potential
$W(t,x)=\omega(t)V(x)$. We apply the theorems to the case of a time dependent
central potential and the harmonic oscillator and determine all Lie and
Noether point symmetries. Finally we prove that these theorems also apply to
the case of a dynamical system with linear dumping and study two examples.

\end{abstract}


PACS - numbers: 2.40.Hw, 4.20.-q, 4.20.Jb, 04.20.Me, 03.20.+i, 02.40.Ky

\section{Introduction}

\label{Introduction}

Newtonian dynamical systems share two fundamental characteristics; the
equation of motion follows from Newton's Second Law and its particular form
depends upon the applied force or potential and the motion occurs in Euclidian
three-dimensional space. This scenario goes over to Special Relativity (for
pure forces)\ by the generalization of Newton's Second Law to Minkowski space.
\ In the case of General Relativity Newton's Law is replaced by the Principle
of Inertia according to which free fall takes place along geodesics. We note
that in all these nonquantum theories we have two common elements:\ A
second-order equation of motion and a (finite-dimensional) Riemannian manifold
in which the motion is considered.

This remark takes us to the next step. The equation of motion as a
second-order differential equation (ODE) admits certain Lie/Noether \ point
symmetries which characterize this differential equation \cite{sLie,stephani}.
On the other hand the geometry (metric) of the space is fixed (to a certain
extent) by its collineations (Killing vectors (KV), Homothetic vectors (HV),
Conformal Killing vectors (CKV), Affine collineations (AC) and Projective
collineations (PC)) \cite{kri}. Hence what follows is to see if there exists a
relation between the point symmetries of the equations of motion of a
particular system and the collineations which define the underlying geometry.

This has been studied in the recent literature where one finds works which
relate the collineations of the metric with the Lie and the Noether point
symmetries of a general (usually conservative) dynamical system with
second-order differential equations
\cite{Katz74,Katz76,Aminova1995,Aminova2000,FMQ,Bok06,Bok06b,camci1,camci2}
while for some systems of partial differential equations in a Riemannian
manifold the connection between Lie/Noether symmetries and collineations was
the subject of study of \cite{ref1,ref2,ref3,ref4}.

In the present paper we follow on the spirit of our previous works
\cite{paper2,paper3,paper4} and we extend our results to the case of a
dynamical system which moves under the action of a certain type of
time-dependent potential. We derive two theorems which allow the direct
explicit calculation of the Lie and the Noether point symmetries of the
equations of motion of the dynamical system in terms of the collineations of
the metric. We apply these theorems to the well-known cases of the
time-dependent harmonic oscillator \cite{Ray1980,Leach1980,Prince1980}, and
the time-dependent Kepler potential \cite{Prince81}, both in Euclidean space.
As new applications we consider the general central time-dependent potential
and determine the Lie and the Noether point symmetries for motion in Euclidean
space. Finally we consider the harmonic oscillator in a space of constant
curvature and determine the Lie and the Noether point symmetries.

In Section \ref{The Lie symmetries of a time depended dynamical System} we
prove Theorem \ref{The Lie for the time dependent system} which gives the Lie
point symmetries in terms of the special PCs of the metric and their
degeneracies together with conditions / constraints which involve the
potential function defining the specific dynamical system. In Section
\ref{The Noether symmetries of a time depended dynamical System} we prove a
second Theorem \ref{Noether symmetries time dependent} which determines the
Noether conditions in terms of the KVs and the HV of the metric and conditions
which involve the potential function. In Section
\ref{The dampted dynamical sysetm} we show the equivalence of the linearly
damped dynamical system with the time-dependent one we consider and give two
theorems which determine the Lie and the Noether point symmetries of a damped
dynamical system in terms of the collineations of the metric and the damping
function. In Section \ref{Examples section} we consider the application of the
general theorems to various well-known and new cases. In Section
\ref{Conclusions} we discuss our conclusions.

\section{Lie point symmetries for a class of nonautonomous systems}

\label{The Lie symmetries of a time depended dynamical System}

We consider a nonautonomous dynamical system with time-dependent potential,
$W(t,x)=$ $\omega\left(  t\right)  V(x),$ the equations of motion of which
are:%
\begin{equation}
\ddot{x}^{i}+\Gamma_{jk}^{i}\dot{x}^{j}\dot{x}^{k}+\omega\left(  t\right)
V^{,i}=0~~,\omega_{,t}V^{,i}\neq0,V(x^{i}). \label{LPT.01}%
\end{equation}

In \cite{paper2}\ the Lie point symmetry conditions for the system of ODEs:
\begin{equation}
\ddot{x}^{i}+\Gamma_{jk}^{i}\dot{x}^{j}\dot{x}^{k}+\sum\limits_{m=0}%
^{n}P_{j_{1}...j_{m}}^{i}\dot{x}^{j_{1}}\ldots\dot{x}^{j_{m}}=0, \label{de.1}%
\end{equation}
where $\Gamma_{jk}^{i}$ are the connection coefficients in an affine space and
$P_{j_{1}...j_{m}}^{i}(t,x)$ are smooth polynomials completely symmetric in
all lower indices have been derived.

Comparing (\ref{LPT.01}) and (\ref{de.1}) we find that the results of the
general case apply for $P^{i}\neq0$ and the remaining $P$- terms zero. We
infer that the Lie point symmetry conditions for a general $P^{i}$ are the following:%

\begin{equation}
L_{\eta}P^{i}+\xi P_{~,t}^{i}+2\xi,_{t}P^{i}+\eta^{i},_{tt}=0 \label{LTP.03}%
\end{equation}%
\begin{equation}
\left(  \xi,_{k}\delta_{j}^{i}+2\xi,_{j}\delta_{k}^{i}\right)  P^{k}+2\eta
^{i},_{t|j}-\xi,_{tt}\delta_{j}^{i}=0 \label{LTP.04}%
\end{equation}%
\begin{equation}
L_{\eta}\Gamma_{(jk)}^{i}=2\xi,_{t(j}\delta_{k)}^{i} \label{LTP.05}%
\end{equation}%
\begin{equation}
\xi_{(,i|j}\delta_{r)}^{k}=0, \label{LTP.06}%
\end{equation}
where $X=\xi(t,x^{j})\partial_{t}+\eta^{i}(t,x^{j})\partial_{x^{i}}$ is the
Lie point symmetry vector. For the specific case we are considering $P\left(
t,x^{k}\right)  ^{i}=\omega\left(  t\right)  V^{,i}\left(  x^{k}\right)  $ the
system of the Lie point symmetry conditions for the equation of motion
(\ref{LPT.01}) reduces to the following:%

\begin{equation}
\omega L_{\eta}V^{,i}+\xi\omega_{,t}V^{,i}+2\omega\xi,_{t}V^{,i}+\eta
^{i},_{tt}=0 \label{LTP1.09}%
\end{equation}%
\begin{equation}
\omega\left(  \xi,_{k}\delta_{j}^{i}+2\xi,_{j}\delta_{k}^{i}\right)
V^{,k}+2\eta^{i},_{t|j}-\xi,_{tt}\delta_{j}^{i}=0 \label{LTP1.10}%
\end{equation}%
\begin{equation}
L_{\eta}\Gamma_{(jk)}^{i}=2\xi,_{t(j}\delta_{k)}^{i} \label{LTP1.11}%
\end{equation}

\begin{equation}
\xi_{(,i|j}\delta_{r)}^{k}=0. \label{LTP1.12}%
\end{equation}
We omit the solution of the determining equations. \ The steps that we follow
are those described in \cite{paper2} for the case of the autonomous system.

\begin{theorem}
\label{The Lie for the time dependent system}The Lie point symmetries of the
time-dependent dynamical system:%
\begin{equation}
\ddot{x}^{i}+\Gamma_{jk}^{i}\dot{x}^{j}\dot{x}^{k}+\omega\left(  t\right)
V^{,i}=0~~~,~\omega_{,t}\neq0, \label{PL.1}%
\end{equation}
are generated by the special PCs of the metric $g_{ij}$ of the Riemannian
space in which the motion takes place, as follows:

\textbf{Case I}

The Lie point symmetries generated by the Affine algebra $Y^{i}$ of the
metric. There are two cases to consider:

\qquad\textbf{Case I.1} The time function $\omega\left(  t\right)  =\frac
{1}{\left(  d_{1}t+d_{2}\right)  ^{2}}.$ The Lie point symmetry for any
potential $V(x^{i})$ is:%
\begin{equation}
X=\left(  d_{1}t+d_{2}\right)  \partial_{t} \label{PL.10}%
\end{equation}

\qquad\textbf{Case I.2} The Lie point symmetry for $\omega\left(  t\right)
\neq\frac{1}{\left(  d_{1}t+d_{2}\right)  ^{2}}$is of the form:%
\begin{equation}
X=\left(  d_{1}t+d_{2}\right)  \partial_{t}+a_{1}Y^{i}, \label{PL.11}%
\end{equation}
where the constants $a_{1},d_{1}$ are computed from the condition:%
\[
\left(  d_{1}t+d_{2}\right)  (\ln\omega)_{,t}=d_{1}-2d_{2}%
\]

and the potential satisfies the condition:%
\begin{equation}
a_{1}L_{Y}V^{,i}+d_{1}V^{,i}=0. \label{PL.9}%
\end{equation}

\textbf{Case II}

The Lie point symmetries are generated by the gradient KVs and/or gradient HV
$Y^{i}$of the metric and $Y^{i}\neq V^{,i}.$

The Lie symmetry vector is:%
\begin{equation}
X=D\left(  t\right)  \partial_{t}+T(t)Y^{i}\partial_{i}, \label{PL.2}%
\end{equation}

where the functions $D\left(  t\right)  ,T\left(  t\right)  $ are computed
form the equations:%
\begin{align}
D(\ln\omega)_{,t}+2D_{,t}  &  =d_{0}T\label{PL.3}\\
T\,_{,tt}  &  =m\omega T \label{PL.4}%
\end{align}
\newline and the potential $V^{,i}$ satisfies the condition:%
\begin{equation}
L_{Y}V^{,i}+d_{0}V^{,i}+mY^{i}=0. \label{PL.5}%
\end{equation}
\textbf{Case III}.

The Lie point symmetries are generated by the vector $Y^{i}=kV^{,i}$ which is
a gradient KV or a gradient HV of the metric$.$

The Lie symmetry vector is:%
\begin{equation}
X=D\left(  t\right)  \partial_{t}+T(t)Y^{i}\partial_{i} \label{PL.6}%
\end{equation}

where the functions $D\left(  t\right)  ,T\left(  t\right)  $ are computed
form the equations:%
\begin{align}
D(\ln\omega)_{,t}+2D_{,t}+\frac{k}{\omega}T  &  =0\label{PL.7}\\
D\,_{,tt}  &  =2\psi T_{,t} \label{PL.8}%
\end{align}
and $\psi=0$ for a KV and $\psi=1$ for a HV.

\textbf{Case IV }

The Lie point symmetries are generated from the proper special projective
vector $Y^{i}.$

The Lie symmetry vector is:%
\begin{equation}
X=\left(  C(t)S_{J}+D(t)\right)  \partial_{t}+T(t)S_{J}V^{,i}\partial_{i}
\label{PL.12}%
\end{equation}
where the functions $C\left(  t\right)  ,D\left(  t\right)  ,T\left(
t\right)  $ satisfy the following system of equations:%
\begin{align}
D\left(  \ln\omega\right)  _{,t}+\frac{2D_{,t}}{T}  &  =0\label{PL.14}\\
D_{,tt}  &  =0\label{PL.15}\\
\frac{C}{T}\left(  \ln\omega\right)  _{,t}+2\frac{C_{,t}}{T}+\frac{\lambda
}{\omega T}T_{,tt}  &  =\lambda_{1}\label{PL.17}\\
T_{,t}  &  =a_{7}\omega C\label{PL.20}\\
C_{,t}  &  =a_{0}T \label{PL.21}%
\end{align}
and the potential $V^{,i}$ satisfies the following conditions:%
\[
\text{ }V^{i}\text{ is a gradient HV}%
\]%
\begin{align}
S_{J,j}V^{,j}-\lambda_{1}S_{J}  &  =0\label{PL.18}\\
\left(  S_{J,k}\delta_{j}^{i}+2S_{J},_{j}\delta_{k}^{i}\right)  V^{,k}%
+a_{7}\left(  2Y_{J}^{;i}{}_{;~j}-a_{0}S_{J}\delta_{j}^{i}\right)   &  =0
\label{PL.19}%
\end{align}
and $S_{J}^{,i}$ are the gradient KVs of the metric.
\end{theorem}

We continue our analysis with the Noether point symmetries.

\section{Noether point symmetries for a class of nonautonomous systems}

\label{The Noether symmetries of a time depended dynamical System}

Consider a particle moving in a space with metric $g_{ij}$ under the influence
of the potential $V\left(  t,x^{k}\right)  .$ Then the Lagrangian describing
the motion of the particle is:
\begin{equation}
L=\frac{1}{2}g_{ij}\dot{x}^{i}\dot{x}^{j}-V\left(  t,x^{k}\right)  .
\label{TP.1}%
\end{equation}
A vector field in the space, $\mathbf{X}=\xi\left(  t,x^{k}\right)
\partial_{t}+\eta^{i}\left(  t,x^{k}\right)  \partial_{x^{i}},$ is a Noether
point symmetry of the Lagrangian if the following condition is satisfied:%
\begin{equation}
\mathbf{X}^{\left[  1\right]  }L+\frac{d\xi}{dt}L=\frac{df}{dt}, \label{TP.1a}%
\end{equation}
where $\mathbf{X}^{\left[  1\right]  }=\xi\left(  t,x^{k}\right)  \partial
_{t}+\eta^{i}\left(  t,x^{k}\right)  \partial_{x^{i}}+\left(  \frac{d\eta^{i}%
}{dt}-\dot{x}^{i}\frac{d\xi}{dt}\right)  \partial_{\dot{x}^{i}}$ is the first
prolongation of $\mathbf{X}.$ We compute:%
\begin{align}
X^{\left[  1\right]  }L  &  =\left(  \xi\partial_{t}+\eta^{i}\partial_{x^{i}%
}+\left(  \frac{d\eta^{k}}{dt}-\dot{x}^{k}\frac{d\xi}{dt}\right)
\partial_{\dot{x}^{k}}\right)  \left(  \frac{1}{2}g_{ij}\dot{x}^{i}\dot{x}%
^{j}-V\left(  t,x^{k}\right)  \right) \nonumber\\
&  =\frac{1}{2}\left(
\begin{array}
[c]{c}%
\eta^{k}g_{ij,k}\dot{x}^{i}\dot{x}^{j}+2\frac{\partial\eta^{i}}{\partial
t}g_{ij}\dot{x}^{j}+\frac{\partial\eta^{i}}{\partial x^{r}}g_{ik}\dot{x}%
^{k}\dot{x}^{r}+\\
+\frac{\partial\eta^{j}}{\partial x^{r}}g_{kj}\dot{x}^{k}\dot{x}^{r}%
-2g_{ij}\dot{x}^{i}\dot{x}^{j}\frac{\partial\xi}{\partial t}-2\frac
{\partial\xi}{\partial x^{k}}g_{ij}\dot{x}^{i}\dot{x}^{j}\dot{x}^{k}%
-2V_{,k}\eta^{k}-2\xi V_{,t}%
\end{array}
\right)  .
\end{align}
The second term gives:%
\begin{equation}
\left(  \xi_{,t}+\dot{x}^{r}\xi_{,r}\right)  \left(  \frac{1}{2}g_{ij}\dot
{x}^{i}\dot{x}^{j}-V\right)  =\frac{1}{2}\left(  \xi_{,t}g_{ij}\dot{x}^{i}%
\dot{x}^{j}-2V\xi_{,t}+g_{ij}\xi_{,k}\dot{x}^{i}\dot{x}^{j}\dot{x}^{k}%
-2\dot{x}^{k}\xi_{,k}V\right)
\end{equation}
and similarly the third term:%
\begin{equation}
\frac{df}{dt}=f_{t}+\dot{x}^{k}f_{,k}.
\end{equation}
When we substitute into (\ref{TP.1a}) and note that the resulting equation is
an identity in $\dot{x}^{k},$ we set the coefficient of each power of $\dot
{x}^{k}$ equal to zero and find the following set of equations:
\begin{align}
V_{,k}\eta^{k}+V\xi_{,t}+\xi V_{,t}  &  =-f_{,t}\label{TP.1b}\\
\frac{\partial\eta^{i}}{\partial t}g_{ij}-\xi_{,j}V  &  =f_{,j}\label{TP.1c}\\
L_{\eta}g_{ij}  &  =2\left(  \frac{1}{2}\frac{\partial\xi}{\partial t}\right)
g_{ij}\label{TP.1d}\\
\frac{\partial\xi}{\partial x^{k}}  &  =0. \label{TP.1e}%
\end{align}

Equation (\ref{TP.1d}) implies $\xi=\xi\left(  t\right)  $ and the system of
the equations reduces as follows:%
\begin{align}
L_{\eta}g_{ij}  &  =2\left(  \frac{1}{2}\xi_{,t}\right)  g_{ij}\label{TP.2}\\
V_{,k}\eta^{k}+V\xi_{,t}+\xi V_{,t}  &  =-f_{,t}\label{TP.3}\\
\eta_{i,t}  &  =f_{,i}. \label{TP.4}%
\end{align}

The system of equations (\ref{TP.2}) - (\ref{TP.4}) for potential functions of
the form:
\begin{equation}
V\left(  t,x^{k}\right)  =\omega\left(  t\right)  V\left(  x^{k}\right)  ~
\label{TP.5}%
\end{equation}
becomes:%
\begin{align}
L_{\eta}g_{ij}  &  =2\left(  \frac{1}{2}\xi_{,t}\right)  g_{ij}\label{TP.6}\\
V_{,k}\eta^{k}+\left(  \xi_{,t}+\left(  \ln\omega\right)  _{,t}\xi\right)  V
&  =-\frac{f_{,t}}{\omega}\label{TP.7}\\
\eta_{i,t}  &  =f_{,i}. \label{TP.8}%
\end{align}
The solution of the system of equations (\ref{TP.6})-(\ref{TP.8}) which leads
to the following theorem.

\begin{theorem}
\label{Noether symmetries time dependent} For a dynamical system with
Lagrangian:
\begin{equation}
L=\frac{1}{2}g_{ij}\dot{x}^{i}\dot{x}^{j}-\omega\left(  t\right)  V\left(
x^{k}\right)  ~,~\omega_{,t}\neq0, \label{Pr.1}%
\end{equation}
the Noether point symmetries are generated by the KVs and the HV of the metric
as follows:\newline\textbf{Case I} \newline Noether point symmetries generated
by the KVs and the HV $Y^{i}.$

The Noether point symmetry vector and the Noether function are:%
\begin{equation}
X=\left(  2a_{1}\psi_{Y}t+d_{1}\right)  \partial_{t}+a_{1}Y^{i}\partial
_{i}~,~~f=c_{2}\int\omega dt, \label{Pr.2}%
\end{equation}
where the constants are computed form the relation:%
\begin{equation}
\left(  \ln\omega\right)  _{,t}\left(  2\psi_{Y}t+d_{2}\right)  =d_{1}
\label{Pr.3b}%
\end{equation}
and the potential function $V(x^{i})$ satisfies the condition:%
\begin{equation}
V_{,k}Y^{k}+d_{2}V+c_{2}=0. \label{Pr.3}%
\end{equation}
$\psi_{Y}$ is a constant which equals zero for a KV and nonzero constant for a
HV.\newline\newline\textbf{Case II} \newline Noether point symmetries produced
by gradient KVs and/or the gradient HV $Y^{i}:$

We have two cases:

\qquad\textbf{Case II.a} $V^{,i}$ is not a gradient HV

The Noether point symmetry vector and the Noether function are:%
\begin{equation}
X=\xi(t)\partial_{t}+T(t)S_{J}^{i}\partial_{i}~,~f=T_{,t}S, \label{Pr.5}%
\end{equation}
where the index $J$ runs through all the gradient KVs. The functions
$\xi(t),T(t)$ and the constants are computed from the equations:%
\begin{align}
\xi_{,t}  &  =2\psi_{Y}T\label{Pr.7}\\
\frac{1}{\omega}\frac{T_{,tt}}{T}  &  =m\label{Pr.8}\\
\frac{1}{T}\left(  \ln\omega\right)  _{,t}\xi &  =d_{1}\label{Pr.9}\\
\frac{1}{T}\frac{K_{,t}}{\omega}  &  =k
\end{align}
and the potential function satisfies the condition:
\begin{equation}
V_{,k}Y^{k}+2\psi_{Y}V+d_{1}V+mS+k=0. \label{Pr.6}%
\end{equation}
\newline\qquad\textbf{Case II.b} \newline$V^{,i}$ is a gradient KV

The Noether point symmetry vector and the Noether function are:
\begin{equation}
X=\xi(t)\partial_{t}+T(t)V^{i}\partial_{i}~,~f=T_{,t}S+K\left(  t\right)  ,
\label{Pr.10}%
\end{equation}
where the functions $\xi(t),T(t)$ are computed from the equations:%
\begin{align}
\xi_{,t}  &  =2\psi_{Y}T\label{Pr.12}\\
\lambda\frac{1}{T}\frac{T_{,tt}}{\omega}+\frac{\left(  \ln\omega\right)
_{,t}\xi}{T}  &  =d_{2}\label{Pr.13}\\
\frac{1}{T}\frac{K_{,t}}{\omega}  &  =k \label{Pr.14}%
\end{align}
and the potential function satisfies the condition:
\begin{equation}
V_{,k}Y^{k}+(2\psi_{Y}+d_{2})V+k=0. \label{Pr.15}%
\end{equation}

\end{theorem}

\section{The equivalence of a linearly damped dynamical system and a
time-dependent dynamical system}

\label{The dampted dynamical sysetm}

Consider a particle moving in a Riemannian space with metric function $g_{ij}$
under the combined action of a linear damping force, $\phi\left(  t\right)
\dot{x}^{i},$ and a conservative force with potential function, $V\left(
x^{i}\right)  .$ The equation of motion of the particle is:%
\begin{equation}
\frac{d^{2}x^{i}}{dt^{2}}+\Gamma_{jk}^{i}\frac{dx^{j}}{dt}\frac{dx^{k}}%
{dt}+\phi\left(  t\right)  \frac{dx^{i}}{dt}+V^{,i}=0. \label{AP.08}%
\end{equation}
If we consider the change of parameter \thinspace$t\rightarrow s=S(t),$ the
equation of motion is written as follows:
\begin{equation}
\frac{d^{2}x^{i}}{ds^{2}}+\Gamma_{jk}^{i}\frac{dx^{j}}{ds}\frac{dx^{k}}%
{ds}+\frac{1}{\left(  \frac{dS\left(  t\right)  }{dt}\right)  ^{2}}\left(
\left(  \frac{d^{2}S\left(  t\right)  }{dt^{2}}\right)  +\phi\left(  t\right)
\frac{dS\left(  t\right)  }{dt}\right)  \frac{dx^{i}}{ds}+\frac{1}{\left(
\frac{dS\left(  t\right)  }{dt}\right)  ^{2}}V^{,i}=0. \label{AP.09}%
\end{equation}
We define the function $S(t)$ by the requirement:%
\begin{equation}
\left(  \frac{d^{2}S\left(  t\right)  }{dt^{2}}\right)  +\phi\left(  t\right)
\frac{dS\left(  t\right)  }{dt}=0 \label{AP.10}%
\end{equation}
and the equation of motion becomes:
\begin{equation}
\frac{d^{2}x^{i}}{ds^{2}}+\Gamma_{jk}^{i}\frac{dx^{j}}{ds}\frac{dx^{k}}%
{ds}+\frac{1}{\left(  \frac{dS\left(  t\right)  }{dt}\right)  ^{2}}V^{,i}=0.
\label{AP.11}%
\end{equation}
The general solution of (\ref{AP.10}) \ is $S(t)=\int e^{-\int\phi(t)dt}%
dt.~$Because the function $S(t)$ is a point transformation, it must have an
inverse so that $t=S^{-1}\left(  s\right)  ~~,~S\left(  S^{-1}\left(
s\right)  \right)  =s.~$Hence the equation of motion is written as follows:%
\begin{equation}
\frac{d^{2}x^{i}}{ds^{2}}+\Gamma_{jk}^{i}\frac{dx^{j}}{ds}\frac{dx^{k}}%
{ds}+\left(  \frac{dS^{-1}\left(  s\right)  }{ds}\right)  ^{2}V^{,i}=0.
\label{AP.14}%
\end{equation}

We define the function $~$%
\begin{equation}
\omega\left(  s\right)  =\left(  \frac{dS^{-1}\left(  s\right)  }{ds}\right)
^{2} \label{AP.14a}%
\end{equation}
and the equation of motion takes the form:%
\begin{equation}
\frac{d^{2}x^{i}}{ds^{2}}+\Gamma_{jk}^{i}\frac{dx^{j}}{ds}\frac{dx^{k}}%
{ds}+\omega\left(  s\right)  V^{,i}=0, \label{AP.15}%
\end{equation}
which is equation (\ref{LPT.01}) we studied in the previous sections.
Therefore by applying the inverse transformation we are able to compute all
Lie and Noether point symmetries of the linear damped time-dependent dynamical
systems with equation of motion (\ref{AP.08}).

One application of the above general result is the time-dependent linearly
damped harmonic oscillator in a space of constant curvature. The equation of
motion of the time-dependent linear system is:
\begin{equation}
\frac{d^{2}x^{i}}{dt^{2}}+\omega\left(  t\right)  x^{i}=0~,~\omega\left(
t\right)  =\frac{\gamma^{2}}{t^{2}}. \label{AP.15a}%
\end{equation}
In this case we have from (\ref{AP.14a}) the transformation:%
\begin{equation}
\left(  \frac{dS^{-1}\left(  t\right)  }{dt}\right)  =\frac{\gamma}%
{t}\rightarrow S^{-1}\left(  t\right)  =\gamma\ln t=s.
\end{equation}
The function $S(t)$ is $1:1$ and the inverse is:%
\begin{equation}
t=e^{\frac{1}{\gamma}s}.
\end{equation}
Therefore:
\[
\frac{d}{dt}\left(  \frac{dx^{i}}{ds}\frac{ds}{dt}\right)  =\frac{d}%
{dt}\left(  \frac{dx^{i}}{ds}\right)  \left(  \gamma\frac{d\left(  \ln
t\right)  }{dt}\right)  ^{2}+\frac{dx^{i}}{ds}\left(  \gamma\frac{d^{2}%
}{dt^{2}}\left(  \ln t\right)  \right)  .
\]
Hence equation (\ref{AP.15a}) becomes:
\[
\frac{d^{2}x^{i}}{ds^{2}}\frac{\gamma^{2}}{t^{2}}-\frac{dx^{i}}{ds}%
\frac{\gamma}{t^{2}}+\frac{\gamma^{2}}{t^{2}}x=0,
\]
that is:
\begin{equation}
\frac{d^{2}x^{i}}{ds^{2}}-\frac{1}{\gamma}\frac{dx^{i}}{ds}+x^{i}=0
\label{Ex.01}%
\end{equation}
which is the damped time-dependent harmonic oscillator. We conclude that the
time-dependent harmonic oscillator has the same Lie and Noether point
symmetries with the linear damped time-dependent harmonic
oscillator~\cite{Mah93,Cervero}.

\section{Applications}

\label{Examples section}

In order to demonstrate our results in this section we determine the point
symmetries of some well-known systems.

\subsection{The Lie and the Noether point symmetries of a generalized
time-dependent central potential}

We determine the Lie and the Noether point symmetries of the time-dependent
central potential $V(t,r)=\omega\left(  t\right)  \frac{1}{n}r^{n}$ in
Newtonian space. The equation of motion is:%
\begin{equation}
\ddot{x}^{i}+\omega\left(  t\right)  \frac{1}{n}r^{n}x^{i}=0~,~\left(
n\neq0,-2,2\right)  . \label{GCP.1}%
\end{equation}
The case $n=2$ corresponds to the harmonic oscillator and the case $n=-2$ is a
singular case. The Kepler potential follows for the value $n=1$.

\subsubsection{Lie point symmetries}

We apply Theorem \ref{The Lie for the time dependent system} in which only
Case II is applied. Condition (\ref{PL.3}) is satisfied by the nongradient KVs
$X_{IJ}~$for $d_{1}=0~$and by the gradient HV $~H^{i},~d_{1}^{\prime}~=\left(
2-n\right)  .$ Therefore we have the following Lie point symmetries:%
\[
X=a_{1}\left(  d_{2}t+d_{3}\right)  \partial_{t}+\left(  a_{0}X_{IJ}%
+a_{1}H^{i}\right)  \partial_{i},
\]
where the constants $d_{2},d_{3}$ must satisfy the condition:%
\[
\left(  d_{2}t+d_{3}\right)  \left(  \ln\omega\right)  _{,t}=\left(
2-n\right)  -2d_{2}.
\]
This condition is satisfied only for $\omega\left(  t\right)  =t^{a},~$where
$a\neq-2,~~$that is, for $\left(  \ln\omega\right)  _{,t}=\frac{a}{t}.$ In
this case $d_{3}=0$ and $d_{2}=\frac{\left(  2-n\right)  }{a+2}.$ We conclude
that Case B possesses the Lie point symmetries:%
\[
X=a_{1}\left(  \frac{\left(  2-n\right)  }{a+2}t\right)  \partial_{t}+\left(
a_{0}X_{IJ}+a_{1}H^{i}\right)  \partial_{i}%
\]
only when $\omega\left(  t\right)  =t^{a},~$where $a\neq-2.$ It is a standard
calculation to show that the finite point transformation generated by the Lie
point symmetry vector corresponding to the parameter $a_{1}$ leaves invariant
the quantity:%
\begin{equation}
\frac{r^{\left(  2-n\right)  }}{t^{\left(  a+2\right)  }}=\text{constant.}
\label{GCP.2}%
\end{equation}
We conclude that for $\omega\left(  t\right)  =t^{a}~$where $a\neq-2$ Kepler's
Third Law is generalized to the time-dependent central potential,
$V(t,r)=t^{a}r^{n-2},~\left(  a\neq-2,n\neq0,-2,2\right)  $.

\subsubsection{Noether point symmetries}

For the Lagrangian function,%
\begin{equation}
L=\frac{1}{2}\delta_{ij}\dot{x}^{i}\dot{x}^{j}-\frac{1}{n}r^{n}, \label{GCP.3}%
\end{equation}
we apply Theorem \ref{Noether symmetries time dependent} and from Case I we
have that condition (\ref{Pr.3}) is met by the nongradient KVs $X_{IJ}~$for
$d_{1}=0~$and by the gradient HV $~H^{i}~$for $~d_{1}~=\left(  n+2\right)  .$
Therefore we have the Noether point symmetries:%
\[
X=a_{1}\left(  2t+d_{2}\right)  \partial_{t}+\left(  a_{0}X_{IJ}+a_{1}%
H^{i}\right)  \partial_{i}%
\]
provided the constant $d_{2}$ satisfies the condition:%
\[
\left(  \ln\omega\right)  _{,t}\left(  2t+d_{2}\right)  =-\left(  n+2\right)
.
\]
This is possible only when $\omega\left(  t\right)  =t^{-\frac{n+2}{2}}$, that
is $\left(  \ln\omega\right)  _{,t}=-\frac{n+2}{2t}$ and $d_{2}=0.$ For this
particular case we have the Noether point symmetries:
\[
X=a_{1}2t\partial_{t}+\left(  a_{0}X_{IJ}+a_{1}H^{i}\right)  \partial_{i}.
\]
We note that the Noether point symmetries are Lie point symmetries when
$a=-\frac{n+2}{2}$ and $n\neq4$.

For the special case, $\omega\left(  t\right)  =t^{-\frac{n+2}{2}},$ the
Noether currents are:%
\begin{align*}
I_{0}  &  =X_{IJ}^{i}\dot{x}_{j}\\
I_{1}  &  =2t\left(  \frac{1}{2}\delta_{ij}\dot{x}^{i}\dot{x}^{j}+\frac{1}%
{n}r^{n}\right)  -H^{i}\dot{x}_{i}.
\end{align*}
The Noether point symmetry corresponding to the parameter $a_{1}$ gives the
first integral (conserved quantity):%
\begin{equation}
\frac{r^{2}}{t}=\text{constant} \label{GCP.4}%
\end{equation}
which follows form (\ref{GCP.2}) if we set $a=-\frac{n+2}{2}.$

\subsection{The Lie and the Noether point symmetries of the exceptional case
$W(r)=\omega\left(  t\right)  r^{-4}$}

This is the one special case we did not study in the last example. The
equation of motion is:
\begin{equation}
\ddot{x}^{i}+\omega\left(  t\right)  x^{i}r^{-4}=0.~ \label{GCP.5}%
\end{equation}

\subsubsection{Lie point symmetries}

Following Theorem \ref{The Lie for the time dependent system} and from Case II
it follows that condition (\ref{PL.3}) is satisfied by the nongradient KVs
$X_{IJ}~$for $d_{1}=0~$and the gradient HV $H^{i}~,~d_{1}^{\prime}~=4.$
Therefore we have the Lie point symmetries:%
\[
X=a_{1}\left(  d_{2}t+d_{3}\right)  \partial_{t}+\left(  a_{0}X_{IJ}%
+a_{1}H^{i}\right)  \partial_{i}%
\]
provided the constants $d_{2},d_{3}$ satisfy the relation:%
\[
\left(  d_{2}t+d_{3}\right)  \left(  \ln\omega\right)  _{,t}=4-2d_{2}.
\]
This is the case only when $\omega\left(  t\right)  =t^{a}~,a\neq-2~~$and
$d_{3}=0,$ $d_{2}=\frac{4}{a+2}.$ In this case we have the Lie point
symmetries:
\begin{equation}
X=a_{1}\left(  \frac{4}{a+2}t\right)  \partial_{t}+\left(  a_{0}X_{IJ}%
+a_{1}H^{i}\right)  \partial_{i}. \label{GCP.6}%
\end{equation}
The Lie point symmetry corresponding to the parameter $a_{1}$ gives the
conserved quantity (generalized Kepler Law):
\begin{equation}
\frac{r^{4}}{t^{\left(  a+2\right)  }}=\text{constant}. \label{GCP.7}%
\end{equation}
In Case III we have the extra Lie point symmetry which is due to the gradient
HV $H^{,i}$ for the values $d_{1}=4~,m=0.$ Then from the solution of the
system of equations (\ref{PL.6}) - (\ref{PL.8}) we compute the functions
$D\left(  t\right)  ,T\left(  t\right)  $ and find the Lie point symmetry
(\ref{PL.5}). In the special case $\omega\left(  t\right)  =t^{a}~,a\neq
-2~~$we compute:%
\[
T\left(  t\right)  =c_{1}~\,,~~D\left(  t\right)  =\frac{4c_{1}}{a+2}t
\]
which is rejected because in this case it is assumed that \ $T_{,t}\neq
0.$\newline We note that the Lie point symmetries of the potential
$V(r)=t^{a}r^{-4}$ $a\neq-2~$follow the same structure with the central
potential $V(r)=t^{a}r^{-n}.$

\subsubsection{\textbf{Noether point symmetries}}

Consider the Lagrangian%
\begin{equation}
L=\frac{1}{2}\delta_{ij}\dot{x}^{i}\dot{x}^{j}+\frac{1}{2}r^{-2}.
\label{GCP.8}%
\end{equation}

From Theorem \ref{Noether symmetries time dependent} for Case I we have that
condition (\ref{Pr.3}) is satisfied by the nongradient KVs $X_{IJ}~$for
$d_{1}=0~$and by the HV $H^{i}~$for$~d_{1}~=0.$ Therefore we have the Noether
point symmetries:%
\[
X=a_{1}\left(  2t+d_{2}\right)  \partial_{t}+\left(  a_{0}X_{IJ}+a_{1}%
H^{i}\right)  \partial_{i}%
\]
provided the constant $d_{2}$ satisfies the relation:%
\[
\left(  \ln\omega\right)  _{,t}\left(  2t+d_{2}\right)  =0
\]
which is impossible. Therefore no Noether symmetry exists in this
case.\newline

For Case II we obtain a Noether point symmetry from the gradient HV $H^{i}%
~$for $d_{1}~=0~,m=0.$ Then the system of equations (\ref{Pr.7}) -
(\ref{Pr.9}) has the trivial solution unless $\omega(t)=$constant. Therefore
we do not have a Noether point symmetry in this case. We conclude that the
Lagrangian (\ref{GCP.8}) has the Noether symmetries:%
\[
X=X_{IJ}%
\]
with corresponding Noether first integrals:%
\[
I_{0}=X_{IJ}^{i}\dot{x}_{j}.
\]
Contrary to the other central potentials in this case we do not have a special
$\omega(t)$ for which a Noether point symmetry is admitted.

\subsection{The time-dependent Kepler potential $W(r)=\omega(t)r^{-2}$}

This is the second singular case to the general central potential. The Lie
point symmetries of the Kepler potential in Euclidian space have been
determined in \cite{Prince81}. The Kepler potential has also been considered
in the case of spaces of constant curvature \cite{Kozlov92} for which the Lie
and the Noether point symmetries have not been determined.

\subsubsection{Lie point symmetries}

The equation of motion is:%
\begin{equation}
\ddot{x}^{i}+x^{i}\frac{\omega\left(  t\right)  }{r^{2}}=0 \label{GCP.9}%
\end{equation}

We apply theorem \ref{The Lie for the time dependent system} and consider
cases. \newline\textit{Case II}

Condition (\ref{PL.3}) is satisfied by the nongradient KVs $X_{IJ}~$for
$d_{1}=0~$and by the HV $~H^{i}~$for $~d_{1}^{\prime}~=3.$ Hence we have the
Lie point symmetries:%
\[
X=a_{1}\left(  d_{2}t+d_{3}\right)  \partial_{t}+\left(  a_{0}X_{IJ}%
+a_{1}H^{i}\right)  \partial_{i},
\]
where the constants $d_{2},d_{3}$ must satisfy the condition:%
\[
\left(  d_{2}t+d_{3}\right)  \left(  \ln\omega\right)  _{,t}=3-2d_{2}.
\]
This is the case only if $\omega\left(  t\right)  =t^{a}~\left(
a\neq-2~\right)  $, $d_{3}=0,d_{2}=\frac{3}{a+2}.$ Therefore we have the Lie
point symmetry:%
\[
X=a_{1}\left(  \frac{3}{a+2}t\right)  \partial_{t}+\left(  a_{0}X_{IJ}%
+a_{1}H^{i}\right)  \partial_{i}%
\]
provided $\omega\left(  t\right)  =t^{a}~\left(  a\neq-2~\right)  .$

One shows easily that under the action of the point transformation generated
by the parameter $a_{1}$ the following quantity is conserved (generalized
Kepler's Third Law):
\begin{equation}
\frac{r^{3}}{t^{\left(  a+2\right)  }}=\text{constant. } \label{GCP.9a}%
\end{equation}

\textit{Case III} and \textit{Case IV} do not give more Lie symmetries.

\subsubsection{Noether point symmetries\newline}

The Lagrangian is:%
\begin{equation}
L=\frac{1}{2}\delta_{ij}\dot{x}^{i}\dot{x}^{j}+r^{-1}. \label{GCP.10}%
\end{equation}

We apply theorem \ref{Noether symmetries time dependent} and consider cases.
\newline\textit{Case I}

Condition (\ref{Pr.3}) is satisfied by the nongradient KVs $X_{IJ}~$for
$d_{1}=0~$and by the HV $~H^{i}~$for $~d_{1}~=q.$ Hence we have the Noether
point symmetries:%
\[
X=a_{1}\left(  2t+d_{2}\right)  \partial_{t}+\left(  a_{0}X_{IJ}+a_{1}%
H^{i}\right)  \partial_{i}%
\]
provided the parameter $d_{2}$ satisfies the condition:%
\[
\left(  \ln\omega\right)  _{,t}\left(  2t+d_{2}\right)  =-1.
\]
This is the case only when $\omega\left(  t\right)  =t^{-\frac{1}{2}}$, that
is $\left(  \ln\omega\right)  _{,t}=-\frac{1}{2t}$and $d_{2}=0.$ The resulting
Noether point symmetry is :
\[
X=a_{1}2t\partial_{t}+\left(  a_{0}X_{IJ}+a_{1}H^{i}\right)  \partial_{i}%
\]
provided $\omega\left(  t\right)  =t^{-\frac{1}{2}}.$ \newline\textit{Case II}
does not give more Noether symmetry vectors.

The First Integrals corresponding to the Noether symmetries are:%
\begin{align}
I_{0}  &  =X_{IJ}^{i}\dot{x}_{j}\label{GCP.11}\\
I_{1}  &  =2t\left(  \frac{1}{2}\delta_{ij}\dot{x}^{i}\dot{x}^{j}%
-r^{-1}\right)  -H^{i}\dot{x}_{i}. \label{GCP.12}%
\end{align}
The Noether symmetry corresponding to the parameter $a_{1}$ gives the
additional conserved quantity (Kepler's Third Law):%
\begin{equation}
\frac{r^{2}}{t}=\text{constant} \label{GCP.13}%
\end{equation}
which is compatible with (\ref{GCP.9a}) when $a=-\frac{1}{2}.$

It is important to note that, when $\omega\left(  t\right)  =t^{\frac{1}{2}},$
Kepler's potential gives the same number of Lie and Noether point symmetries
as the time independent Kepler potential, $\omega\left(  t\right)  =$constant.
Furthermore the Noether point symmetries are the same as the Lie point
symmetries for both cases.

\subsection{The time dependent harmonic oscillator in Euclidean space}

Lastly we consider the time-dependent harmonic oscillator in Euclidean space
which has been extensively studied using the standard Lie approach
\cite{Leach1980},\cite{Prince1980}. Here we use Theorem
\ref{The Lie for the time dependent system} and retrieve these results
directly using the projective algebra of Euclidean space. The equation of
motion is:
\begin{equation}
\ddot{x}^{i}+\omega\left(  t\right)  x=0~,~\text{ }\omega\left(  t\right)
\neq\frac{1}{4t^{2}}. \label{TDHO.1}%
\end{equation}
The projective algebra of Euclidean space is:%
\begin{align*}
S_{I}^{,i}  &  :\text{gradient KV~,~}\mathbf{X}_{IJ}:\text{nongradient KV}\\
H^{,i}  &  :\text{gradient HV~~,~}A_{I}^{i}:\text{Affine Collineation}\\
P_{I}^{i}  &  :\text{special PC.}%
\end{align*}

In Cartesian coordinates these vectors are:%
\begin{align*}
\mathbf{S}_{I}  &  =S_{I}^{,i}\partial_{i}\Rightarrow S_{I}=\delta_{I}%
^{i}\text{ }i=1,2,...,n\\
\mathbf{X}_{IJ}  &  =\delta_{\lbrack I}^{a}\delta_{J]}^{b}x_{a}\partial_{b}\\
\mathbf{H}  &  =x^{i}\partial_{i}~\Rightarrow H\mathbf{=}\frac{1}{2}\left(
x^{i}x_{i}\right) \\
\mathbf{A}_{I}  &  =S_{I}S_{I}^{,i}\partial_{i}=x_{I}\delta_{I}^{i}%
\partial_{i}~\Rightarrow A_{I}^{i}=x_{I}\delta_{I}^{i}\\
\mathbf{P}_{I}  &  =S_{I}\mathbf{H}=x_{I}x^{i}\partial_{i}\Rightarrow
P_{I}^{i}=x_{I}x^{i}.~
\end{align*}
Following Theorem \ref{The Lie for the time dependent system} we consider
cases.\newline\textit{Case II}

Condition (\ref{PL.3}) is satisfied by $X_{IJ}~$\ for $d_{1}=0,$by $H^{,i}$
for $d_{1}^{\prime}=0~$and by $A^{i}~$for$~d_{1}^{^{\prime\prime}}=0$. From
(\ref{PL.4}) follows $D\left(  t\right)  =0.$ Therefore (\ref{PL.2}) gives the
following Lie point symmetry vectors:%
\begin{equation}
X=\left(  l_{1}H^{,i}+l_{2}X_{IJ}+l_{3}A^{i}\right)  \partial_{i}.
\label{TDHO.2}%
\end{equation}
\textit{Case III}

\textit{Case III.1}

We check that the gradient KVs $S_{I}^{,i}$ $J=1,2,\ldots,n$ provides a
symmetry vector for $m=-1~,~d_{0}=0.$We compute $D\left(  t\right)  =0$ and
\ that $T\left(  t\right)  $ is found from the ODE:%
\begin{equation}
T_{,tt}=-\omega(t)T. \label{TDHO.3}%
\end{equation}
The resulting Lie point symmetry vectors are:
\begin{equation}
X=T\left(  t\right)  S_{J}^{,i}\partial_{i}. \label{TDHO.4}%
\end{equation}

\textit{Case III.2}

Because the potential of the harmonic oscillator is the gradient HV we have
the Lie point Symmetry:
\begin{equation}
X=D\left(  t\right)  \partial_{t}+T\left(  t\right)  V^{,i}\partial_{i},
\label{TDHO.5}%
\end{equation}
where the functions $D\left(  t\right)  ,T\left(  t\right)  $ are the solution
of the system of equations:%
\begin{align}
D\left(  \ln\omega\right)  _{,t}+2D_{,t}+\frac{1}{\omega}T_{,tt}  &
=0\label{TDHO.6}\\
2T_{,t}-D_{,tt}  &  =0. \label{TDHO.7}%
\end{align}

\textit{Case IV}.

From equation (\ref{PL.14}) it follows that $D\left(  t\right)  =0$, and
$\lambda=1~,~\lambda_{1}=1~,~\lambda_{2}=-1.$ In this case we also have
$a_{0}=1.$ \ Hence the system of equations (\ref{PL.15})-(\ref{PL.17}) which
yields the functions $C(t),T(t)$ becomes:
\begin{align}
\frac{C}{T}\left(  \ln\omega\right)  _{,t}+\frac{2C_{,t}}{T}+\frac{1}%
{\omega(t)}\frac{T_{,tt}}{T}  &  =1\label{TDHO.8}\\
T_{,t}  &  =-\omega(t)C\label{TDHO.9}\\
C_{,t}  &  =T. \label{TDHO.10}%
\end{align}
The resulting Lie pont symmetry vectors are:%
\[
X=C\left(  t\right)  S_{J}\partial_{t}+T\left(  t\right)  S_{J}V^{,i}\partial
i.
\]

Counting the symmetry vectors we see that they are $N^{2}-\not 1 ,$ $N=n+2$
where $n$~is the dimension of the space. Hence the admitted Lie algebra is the
$sl\left(  n+2,R\right)  $, which means that the system is maximally symmetric
and equivalent with the free particle \cite{Leach1980,Prince1980}.

\section{Conclusions}

\label{Conclusions}

This work extends our previous analysis for the relation between symmetries of
differential equations and collineations of the underlying manifold.
Specifically in \cite{paper2} we studied the connection between projective
collineations and point symmetries of the equations of motion of a particle
moving in a general affine space under the action of an autonomous potential,
we called that $\left(  1\times m\right)  ~$system. Similar analyses have been
done and for linear partial differential equation $\left(  n\times1\right)
~$for which a connection~with the conformal algebra of the space of the
independent variables was found \cite{paper3}. Recently for a system of
quasilinear partial differential equation, $\left(  n\times m\right)  $
systems, and for $m>1$, it was found that that the point symmetries are
related with the conformal algebra of the space of the independent variables
and the affine collineations of the space which defines the evolution of the
dependent variables \cite{paper4}. However, the system of ~$\left(  1\times
m\right)  ~$nonautonomous equations could not be covered from the previous works.

We have proved two general theorems which prove that the Lie and the Noether
point symmetries of a dynamical system moving in Riemannian space under the
action of the time-dependent potential $W(t,x)=\omega(t)V(x)$ are given in
terms of the collineations of the space and a set of constraint conditions
involving the potential and the collineation vectors. We have applied the
theorems to equations of motion in the Euclidian space and derived the Lie and
the Noether point symmetries for the time-dependent central motion and the
time-dependent harmonic oscillator. We have also proved that these theorems
can be used to determine the point symmetries of dynamical systems with linear
dumping. It is interesting to note that this result may be applied to any
second-order equation where the connection coefficients is taken as a general
connection $\Gamma_{jk}^{i}$.

That geometric approach is an alternate method for the study of the symmetries
of differential equations. It is useful because we can study the possible
admitted Lie algebras by just using results from differential geometry and
derive with minimal calculation the collineations of the underlying manifold.
However, except from the technical part it provides a strong relation between
geometry and dynamical systems.

In a forthcoming work we plan to investigate the geometric properties for the
evolution of the dynamical systems by using that geometric point of view.

\bigskip

{\large {\textbf{Acknowledgements}}} \newline The research of AP was supported
by FONDECYT postdoctoral grant no. 3160121. AP thanks the Durban University of
Technology for the hospitality provided while part of this work was performed.


\begin{thebibliography}{99}                                                                                               %


\bibitem {sLie}S. Lie, \textit{Differentialgleichungen}, AMS Chelsea
Publishing, New York \textbf{14} (1967).

\bibitem {stephani}H. Stephani, Differential Equations: Their Solutions using
Symmetry, Cambridge University Press, New York, (1989)

\bibitem {kri}K.L. Duggal and R.Sharma, Symmetries of Spacetimes and
Riemannian Manifolds, Kluwer Academic Press, (1999)

\bibitem {Katz74}G.H. Katzin and J. Levine, \emph{ }J. Math. Phys.\emph{
}\textbf{15,}\emph{\textbf{\ }}1460 (1974)

\bibitem {Katz76}G.H. Katzin and J. Levine, \emph{ }J. Math.
Phys.\textbf{\ 17}, 1345 (1976)

\bibitem {Aminova1995}A.V. Aminova, Sbornik Mathematics \textbf{186} (12),
1711 (1995)

\bibitem {Aminova2000}A.V. Aminova, Tensor, N.S. \textbf{65,} 68, (2000)

\bibitem {FMQ}T. Feroze, F.M. Mahomed and A. Qadir, Nonlinear Dynamics
\textbf{65,} 74 (2006)

\bibitem {Bok06}A.H. Bokhari, A.H. Kara, \ A.R. Kashif and F.D. Zaman, Inter.
Jour. Theor. Phys. \textbf{45}, 1063 (2006)

\bibitem {Bok06b}A.H. Bokhari and A.H. Kara ,Gen Rel. Grav. 39, 2053 (2007)

\bibitem {camci1}U. Camci, Gen. Rel. Grav. 46, 1824 (2014)

\bibitem {camci2}U. Camci, EPJC 74, 3201 (2014)

\bibitem {ref1}Y. Bozhkov and I.L. Freire, J. Differential Equations 249, 872 (2010)

\bibitem {ref2}A. Paliathanasis, M. Tsamparlis and M.T. Mustafa, IJGMMP 12,
1550033 (2015)

\bibitem {ref3}M. Tsamparlis, A. Paliathanasis and A. Qadir, IJGMMP 12,
1550003 (2015)

\bibitem {ref4}S. Jamal, A.H. Kara, A.H. Bokhari, Canadian Journal of Physics,
90, 667 (2012)

\bibitem {paper2}M. Tsamparlis and A. Paliathanasis, Gen. Rel. Grav. 43, 1861 (2011)

\bibitem {paper3}M. Tsamparlis and A. Paliathanasis IJGMMP 11, 1450037 (2014)

\bibitem {paper4}A. Paliathanasis and M. Tsamparlis, J. Geom. Phys. 107, 45 (2016)

\bibitem {Ray1980}J.R. Ray, J. Phys A.: Math. Gen., \textbf{13}, 1969 (1980)

\bibitem {Leach1980}P.G.L. Leach, J Physics A: Math. Gen. \textbf{13}, 1991 (1980)

\bibitem {Prince1980}G.E. Prince and C.J. Eliezer, J. Phys A.: Math. Gen.
\textbf{13} , 815 (1980)

\bibitem {Prince81}G.E. Prince and C.J. Eliezer J. Phys A.: Math. Gen.
\textbf{14}, 587 (1981)

\bibitem {Mah93}F.M. Mahomed, A.H. Kara and P.G.L. Leach, J. Math. Anal.
Appl., \textbf{178}, 116 (1993)

\bibitem {Cervero}J.M. Cervero and J. Villarroel, J. Phys A.: Math. Gen.,
\textbf{ 17}, 1777 (1984)

\bibitem {Kozlov92}V.V. Kozlov and O.A. Harin, Celestial. Mech. Dynam.
Astronom., \textbf{54}, 393 (1992)
\end{thebibliography}
\end{document}